\documentclass[12pt]{amsart}
\usepackage{color}
\usepackage{epsfig}
\usepackage{verbatim}

\usepackage{amssymb}
\usepackage{parskip}

\newtheorem{theorem}{Theorem}[section]

\newtheorem{lemma}[theorem]{Lemma}

\begin{document}

\title {Pancyclicity of almost-planar graphs}
\maketitle

\begin{center}
Santiago T. Adams \footnote{California Institute of Technology, 1200 E California Blvd., Pasadena, CA 91125. Email: stadams@caltech.edu.} and
S. R. Kingan \footnote{Department of Mathematics,
Brooklyn College, 2900 Bedford Avenue, Brooklyn, NY 11210,
and CUNY Graduate Center, 365 Fifth Avenue, New
York, NY 10016. 
Email: skingan@gc.cuny.edu.} 

\end{center}

\begin{abstract} A non-planar graph is almost-planar if either deleting or contracting any edge makes it planar. A graph with $n$ vertices is pancyclic if it contains a cycle of every length from $3$ to $n$, and it is Hamiltonian if it contains a cycle of length $n$. A Hamiltonian path is a path of length $n$ and a graph with a Hamiltonian path between every pair of vertices is called Hamiltonian-connected.  In 1990, Gubser characterized the class of almost-planar graphs. This paper explores the pancyclicity of these graphs. We prove that a $3$-connected almost-planar graph is pancyclic if and only if it has a cycle of length 3. Furthermore, we prove that a 4-connected almost-planar graph is both pancyclic and Hamiltonian-connected.

\end{abstract}

\section {\bf Introduction \label{intro}}

In this paper we examine the pancyclicity and Hamiltonian-connectivity of almost-planar graphs.
A graph on $n$ vertices can have cycles of length $3, 4, \dots ,n$. The \textit{cycle spectrum} of a graph is a list of the cycle lengths present in the graph. A graph that has at least one cycle of every possible length is called \textit{pancyclic}.  A cycle of length $n$  is called a \textit{Hamiltonian cycle} and a path of length $n$ is called a \textit{Hamiltonian path}. A graph is called \textit{Hamiltonian} if it has a Hamiltonian cycle and  \textit{Hamiltonian-connected} if there is a Hamiltonian path between every pair of vertices. Hamiltonian-connectivity implies Hamiltonicity, but not vice-versa. For example, the cycle graph $C_n$ is clearly Hamiltonian, but not Hamiltonian-connected.  

To delete an edge $e$ from a graph $G$ remove it and leave its end vertices intact. The resulting graph, denoted by $G\backslash e$, is called an \textit{edge-deletion} of $G$. To contract an edge $f$ with end vertices $v$ and $v'$, collapse the edge by identifying $v$ and $v'$ as one vertex, and delete the loop formed. The resulting graph, denoted by $G / f$, is called an \textit{edge-contraction}. A graph $H$ is a \textit{minor} of a graph $G$ if $H$ can be obtained from $G$ by deleting or contracting edges.  A graph $H$ is a \textit{deletion-minor} of a graph $G$ if $H$ can be obtained from $G$ by deleting edges.   

A graph is \textit{planar} if it can be drawn in the plane with no pair of edges crossing each other.
Such a drawing is called a \textit{plane representation}.  Otherwise, the graph is \textit{non-planar}. 
 A non-planar graph $G$ is \textit{almost-planar} if  either $G\backslash e$ or $G/e$ is planar for every edge $e$.

A set of vertices whose removal disconnects the graph or results in the trivial graph is called a \textit{vertex-cut}.  A graph with at least $k+1$ vertices is \textit{$k$-connected} if it has no vertex-cut of size less than $k$. In 1956 Tutte proved that a 4-connected planar graph $G$ is Hamiltonian \cite{Tutte1956} and in 1983 Thomassen proved  that $G$ is also Hamiltonian-connected \cite{Thomassen1983}. A 4-connected planar graph is not necessarily pancyclic. For example, the line graph of a cyclically 4-connected planar cubic graph of girth at least 5 is 4-connected and planar, but it does not have a cycle of length 4, as Malkevitch noted in 1972. He conjectured that a 4-connected planar graph is pancyclic if and only if it has a cycle of length 4. Since then, considerable progress on analyzing the cycle spectrum has been made in  for example \cite{ChenFanYu2004},  \cite{Ellingham2019}, \cite{Ozeki2015}, \cite{ThomasYu1994}, and \cite{Thomassen1983}.  In particular, Chen, Fan, and Yu proved that a 4-connected planar graph has a cycle of length $n$, $n-1$, $n-2$, $n-3$, $n-4$, $n-5$, and $n-6$.  \cite{ChenFanYu2004}.  But to-date this fifty-year old conjecture remains unresolved.

Motivated by Tutte's 1956 result and Thomassen's 1983  result that 4-connected planar graphs are Hamiltonian and Hamiltonian-connected, respectively, we set about examining the class of 3-connected almost-planar graphs. The following results are the main results in this paper.

\begin{theorem}  \label{mainresult2}
Let $G$ be a $3$-connected almost-planar graph. Then $G$ is pancyclic if and only if $G$ has a cycle of length $3$.
\end{theorem}

\begin{theorem} \label{mainresult1} Let $G$ be a $4$-connected  almost-planar graph.  Then $G$ is pancyclic and Hamiltonian-connected.
\end{theorem}

The terminology follows \cite{Kingan2022} for the most part.  Section 2 has a description of the 3-connected almost-planar graphs and  Section 3 has all the proofs. The proof of Theorem \ref{mainresult2} is rather long and it  is broken up into several results. Section 4 has concluding remarks and next steps. 


\section{\bf Infinite families of almost-planar graphs \label{almostplanar}}

In 1990 Gubser classified the 3-connected almost-planar graphs into 3 broad classes:  the well-known M{\"o}bius Ladders, and two infinite classes of graphs that we will refer to as the $\mathcal{B}$ graphs and the $\mathcal{H}$ graphs. 

\begin{figure}[h]
\centering
\includegraphics[width=2.9in]{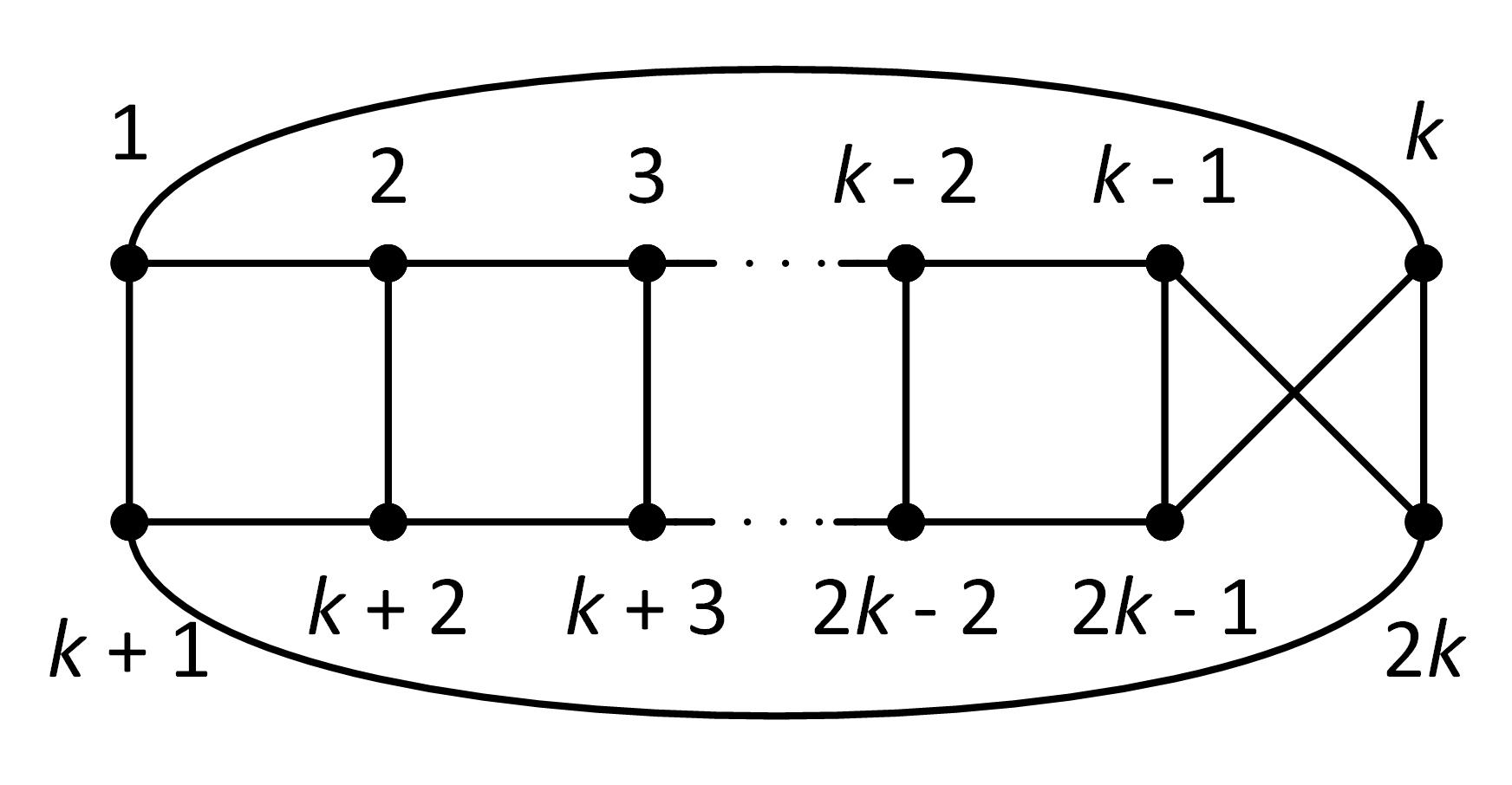}
\caption{The M\"obius ladder $V_{2k}$, where $k\ge 3$.}
 \label{fig1}
\end{figure}

For $k\ge 3$, let  $V_{2k}$ denote the M\"obius ladder on $2k$ vertices   labeled $1, 2, \dots  , k$ along the left and $k+1, \dots ,2k$ along the right, as shown in Figure \ref{fig1}.  Observe that $V_6\cong  K_{3,3}$.

\begin{figure}[h]
\centering
\includegraphics[width=6.2in]{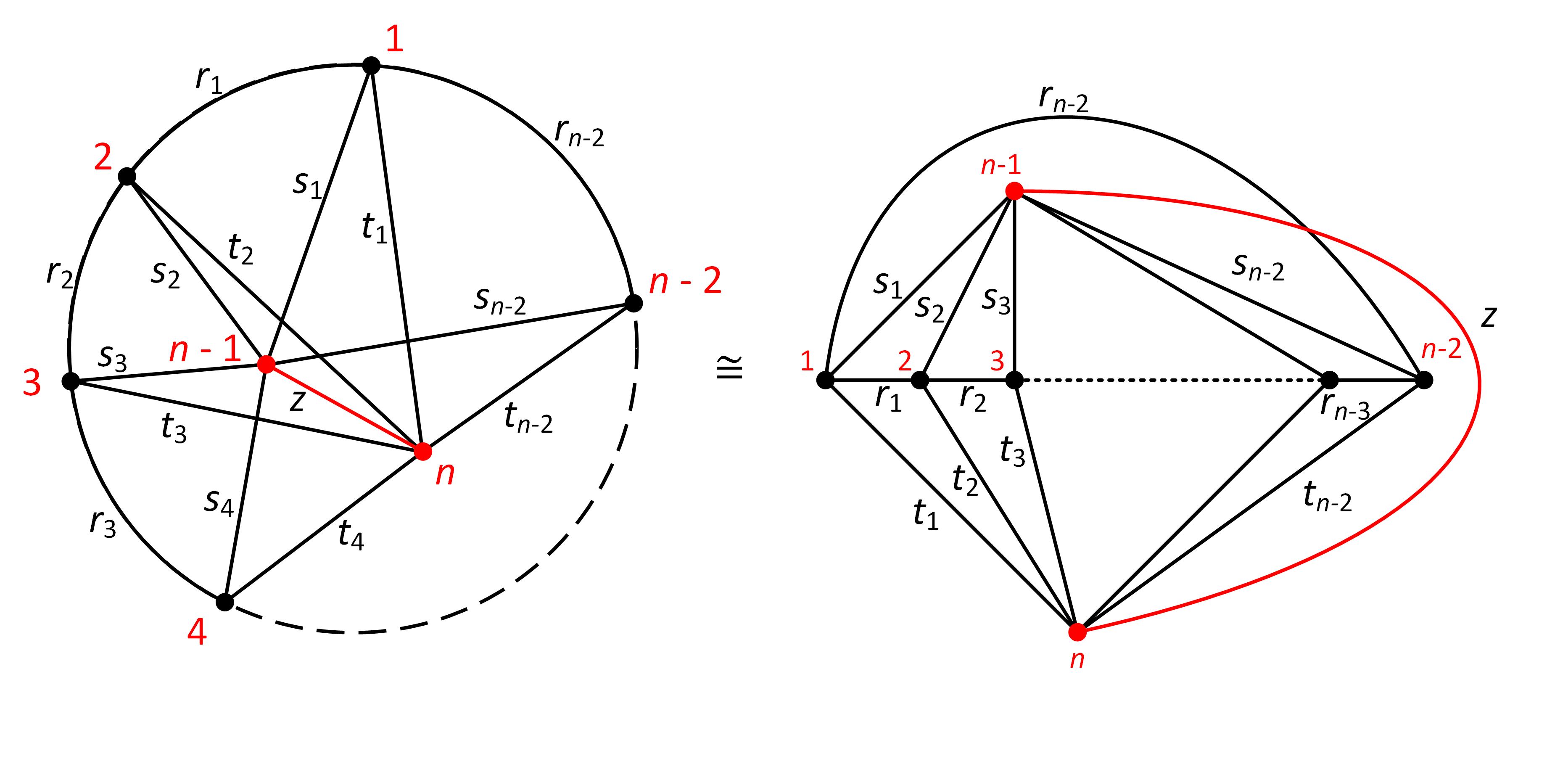}
\caption{The bicycle wheel $B_n$, for $n\ge 5$.}
 \label {fig2}
\end{figure}

For $n\ge 5$, let $B_n$ denote the bicycle wheel with $n$ vertices and $3n-5$ edges shown in Figure  \ref{fig2}.   Label the rim vertices as $1, 2, \dots , n-2$ and the two middle vertices as $n-1$ and $n$.       
Label the $n-2$ edges of the rim as $r_1, \dots , r_{n-2}$.   Label the spokes originating from one of                               
the middle vertices as $s_1, \dots , s_{n-2}$  and from the other as  $t_1, \dots , t_{n-2}$. We will call these spokes \textit{$s$-spokes} and \textit{$t$-spokes}, respectively. Let $z$ be the edge joining the two middle vertices.   An isomorphic drawing of $B_n$, useful for identifying cycles, is also shown in Figure \ref{fig2}.   

\begin{figure}[h]
\centering
\includegraphics[width=5.45in]{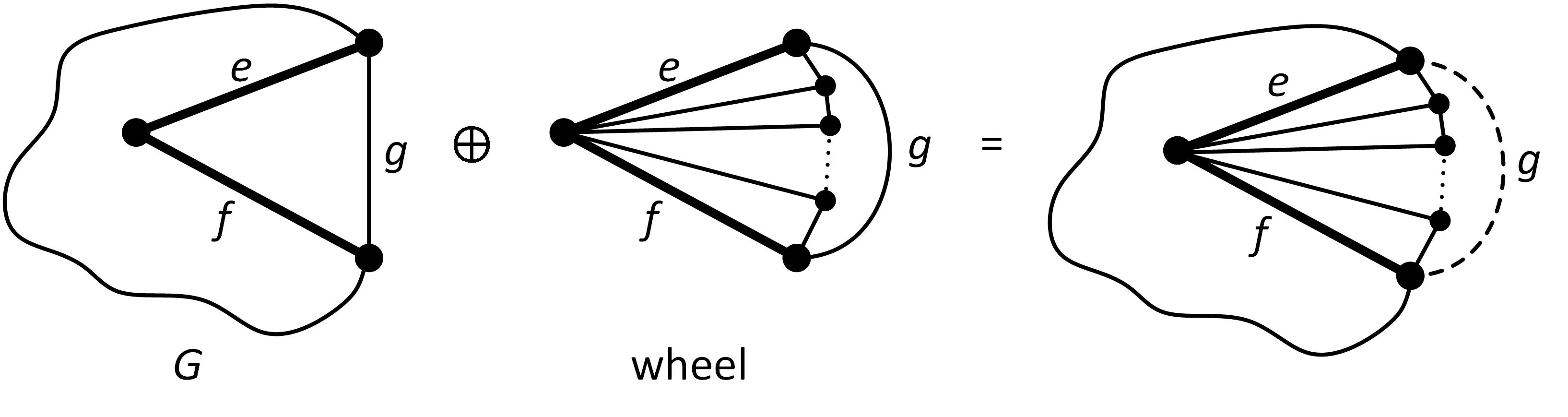}
\caption{The operation of attaching a fan in a triangle}
 \label {fig3}
\end{figure}

A cycle of length $k$ is called a {\it $k$-cycle} in short, and a $3$-cycle is called a {\it triangle}. A {\it clique} in a graph is a maximal induced subgraph isomorphic to the complete graph $K_t$, where $t\ge 2$. Suppose each of the graphs $G_1$ and $G_2$ has a clique $K_t$. The \textit{clique sum} of $G_1$ and $G_2$ across $K_t$ is formed by placing two copies of $K_t$, one on top of the other, and identifying the vertices and edges. For example, let $G$ be a graph with a triangle $K_3$ whose edges are labeled $e, f, g$. Label the edges of any triangle of the wheel graph  as $e, f, g$ so that $e$ and $f$ are spoke edges and $g$ is a rim edge, and consider the clique sum of $G$ and the wheel across the specified triangle as shown in Figure \ref{fig3}. Further delete edge $g$. This operation is called \textit{attaching a type-1 fan} across a triangle in $G$ with \textit{sides} $e$ and $f$. Since this paper talks only about type-1 fans, we will shorten the term to ``fan.''

\begin{figure}[h]
\centering
\includegraphics[width=5.65in]{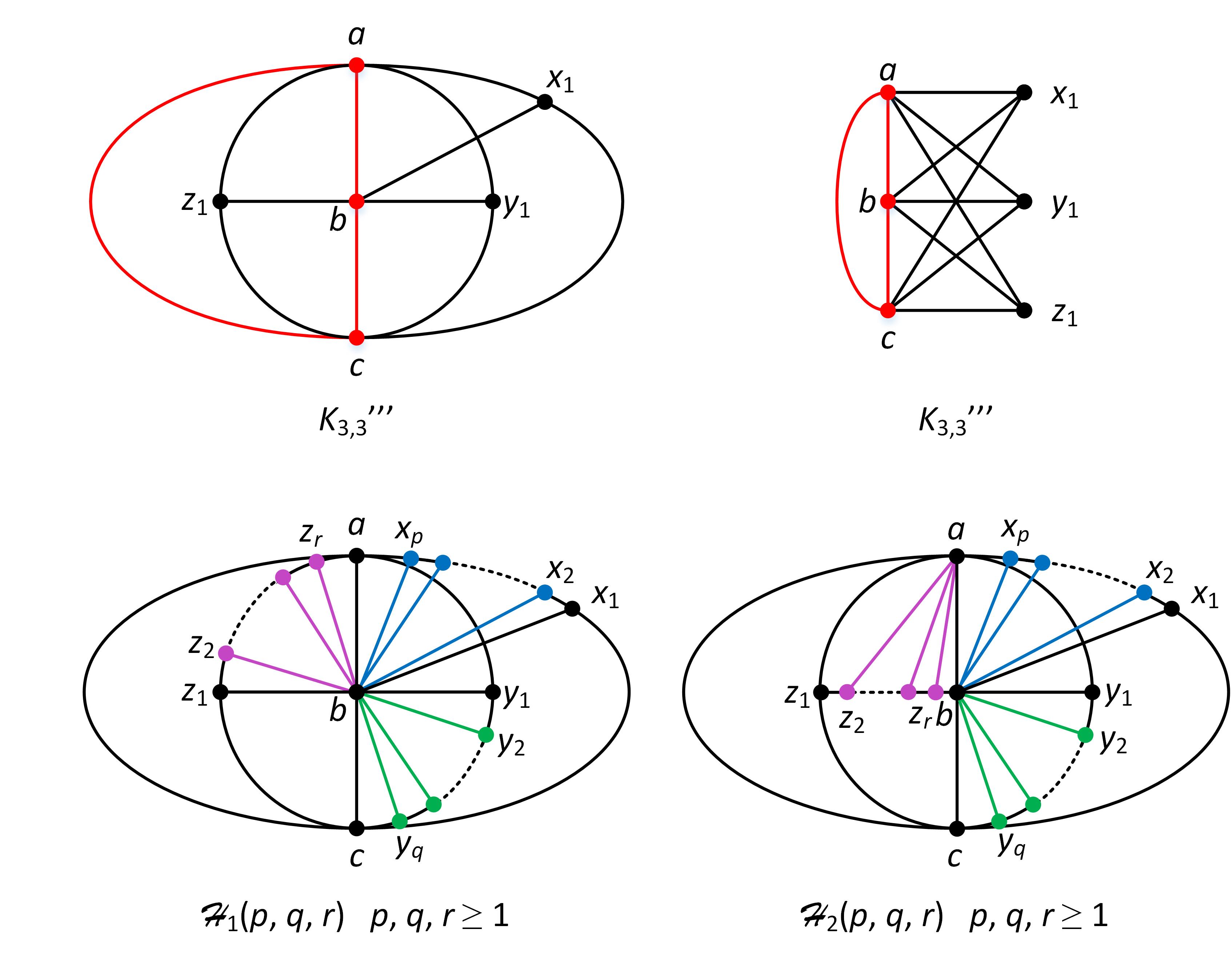}
\caption{Sticking fans in $K_{3,3}'''$}
 \label {fig4}
\end{figure}

Let  $K'''_{3,3}$ be the first graph shown in Figure \ref{fig4} with triangles  $\{a, b, x_1, a\}$,  $\{b, c, y_1, b\}$,   and  $\{a, b,  z_1, a\}$ highlighted.   An isomorphic copy of $K_{3,3}'''$ drawn in the more familiar manner is placed next to it.
 
\begin{itemize}
\item For $p, q, r\ge 1$, let ${\mathcal H_1}(p, q, r)$ be the infinite family  obtained by   attaching  fans of length $p$, $q$, and $r$ along the above  triangles with sides $ab$ and $bx_1$;  $bc$ and  $by_1$; and $ab$ and $\bf bz_1$,  respectively.  
 
\item Let ${\mathcal H_2}(p, q, r)$ be the infinite family obtained by attaching fans of length $p$, $q$, and $r$ along                                                                                                                                                            
 the above triangles with sides $ab$ and $bx_1$;  $bc$ and  $by_1$; and $ab$ and $\bf az_1$,  respectively.     
\end{itemize}
 Note that the only difference in ${\mathcal H_1}(p, q, r)$ and ${\mathcal H_1}(p, q, r)$ is the sides of the fans in triangle  $\{a, b,  z_1, a\}$.  In the first family the sides are $ab$ and $\bf bz_1$ and in the second family the sides are $ab$ and $\bf az_1$.

The following theorem is a restatement of Gubser's result  on  almost-planar graphs \cite{Gubser1996}, taken from \cite{KinganLemos2002}.

\begin{theorem} \label{almostplanartheorem} (Gubser, 1990) A $3$-connected non-planar graph $G$ is almost-planar if and only if $G\cong V_{2k}$, for $k\ge 4$, or $G$ is a $3$-connected deletion-minor of $B_{n}$,  for $n\ge 5$, or of ${\mathcal H_1}(p, q, r)$ or  ${\mathcal H_2}(p, q, r)$, for $p, q, r \ge 1$.  
\end{theorem}

We will refer to the infinite family of bicycle wheels $B_n$ and all its non-planar 3-connected deletion-minors  as the $\mathcal B$ graphs and  the infinite families ${\mathcal H_1}(p, q, r)$ and ${\mathcal H_2}(p, q, r)$, and all their non-planar 3-connected deletion-minors as the  $\mathcal H$ graphs.
 
There exist a number of sufficient conditions for the existence of Hamiltonian cycles based on the degree sequence. The first such result is Dirac's theorem \cite{Dirac1952} which states that if every vertex in a graph $G$ with $n$ vertices has degree at least $\frac{n}{2}$, then $G$ is Hamiltonian.  See also \cite{Ore1960}, \cite{Posa1962}, \cite{Chvatal1972}, \cite{Fan1984}, and \cite{FournierFraisse1985}. In \cite{Bondy1971}, Bondy made the observation that besides  a few  exceptional graphs, degree sequence conditions sufficient for the existence of Hamiltonian cycles are also sufficient for the existence of cycles of all sizes. As a result, there were a number of papers extending  sufficient conditions for Hamiltonicity based on degree sequences to sufficient conditions for pancyclicity, and this has become an active area of research with some recent celebrated breakthroughs. See, for example, \cite{SchmeichelHakimi1974}, \cite{BenhocineWodja1987}, \cite{SchmeichelHakimi1988}, and \cite{Sudakov2023}.
The techniques in these papers  fail for 3-connected almost-planar  graphs since they have at most two vertices whose degrees increase as $n$ increases. Therefore, it is not possible to form a collection of vertices with large enough degree to ensure Hamiltonicity or pancyclicity using existing results. As such, our approach is algorithmic in nature. A description of how to find the cycle together with actually listing the cycle makes it relatively easy for the reader to verify its presence. 

Lastly, we will make frequent use of the well-known result that $G$ is bipartite if and only if $G$ has no cycles of odd length. A cycle of even length is called an {\it even cycle} in short, and a cycle of odd length is called an {\it odd cycle}.

 
\section{\bf Main Results} 

This section has all the new results. Each of the subsections examines the cycle spectrum of an infinite family of almost-planar graphs and all its deletion-minors.  


\subsection{M\"obius Ladders} Each graph in the infinite family of 3-connected graphs $V_{2k}$, where $k\ge 3$, shown in Figure \ref{fig1}, is cubic and therefore has no 3-connected deletion-minors. So we must only analyze the cycle spectrum of $V_{2k}$, which is straightforward. This infinite family has no triangles. Using the twisted ladder representation of $V_{2k}$ illustrated in Figure \ref{fig1}, observe that cycles of length $2t$, where $2\le t \le k$, are obtained by going along one side of the ladder from vertex 1 to vertex $\frac{t}{2}$, crossing a rung of the ladder to vertex $k+\frac{t}{2}$ and returning back to vertex $k+1$ and finally to vertex 1. 
All even cycles can be obtained in this manner. If $k$ is odd, then $V_{2k}$ is a bipartite graph, so it does not have odd cycles. If $k$ is even, then $V_{2k}$ is not a bipartite graph and could have odd cycles, so an additional argument is required.

\begin{theorem} \label{Mobiusgraphs} If $k\ge 4$ is odd, the cycle spectrum of $V_{2k}$  consists of the even numbers from $4$ to $2k$. If $k$ is even, the cycle spectrum of  $V_{2k}$ consists of the even numbers from $4$ to $2k$, as well as all the numbers from $k+1$ to $2k$.
\end{theorem}

\begin{figure}[h]
\centering
\includegraphics[width=3.9in]{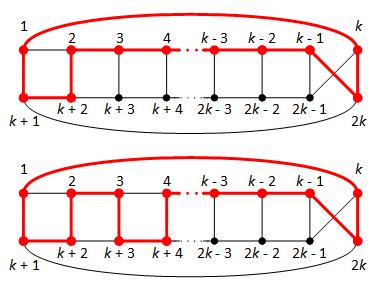}
\caption{Cycles of size $k+3$ and $k+5$ in $V_{2k}$, where $k$ is even}
 \label {V2kEvenCycles}
\end{figure}

\begin{proof}    Observe from Figure \ref{fig1}  that a cycle of length $2t$,  where $2 \le t\le k-1$, is given by
$$1, 2, \dots , \frac{t}{2}, k+\frac{t}{2}, \dots , k+1, 1,$$
and a Hamiltonian cycle is given by
$$1, 2, 3, \dots, k-1, 2k, k, 2k-1, \dots, k+1, 1.$$ 
If $k$ is odd, then $V_{2k}$ is bipartite, and therefore has no odd cycles, so we have a complete classification of its cycle spectrum. If $k$ is even, we will prove that $V_{2k}$ contains all  cycles of length $k+t$, where $t$ is odd and $1\le t \le k-1$. First, observe that since $k$ is even, $k+1, k+3, \dots , k+(k-1)$ are odd. A cycle of length $k+1$ is given by
$$1, 2, 3, \dots , k-1, 2k, k, 1.$$ A cycle of length $k+3$ is given by
$$\ 1, k+1, k+2, 2, 3, 4 , \dots , k-3, k-2, k-1, 2k, k, 1$$
and a cycle of length $k+5$ is given by
$$\ 1, k+1, k+2, 2, 3, k+3, k+4, 4 , \dots , k-3, k-2, k-1, 2k, k, 1.$$ Both are shown in red  in Figure \ref{V2kEvenCycles}. Each additional crossing adds 2 to the length of the cycle. Finally, a  cycle of length $k+(k-1)=2k-1$ is given by
$$1, k+1, k+2, 2, 3, k+3, k+4, 4, 5, 6, \dots , k-3, 2k-3, 2k-2, k-2, k-1, 2k, 1.$$ Lastly, observe from Figure \ref{fig1} that removing edges with end vertices $1$ and $k$ and $k+1$ and $2k$ gives a bipartite graph which has only even cycles. So, any odd cycle in $V_{2k}$ must contain one of the two edges. However, any such cycle has length at least $k$. So there are no odd cycles of length less than $k$. \end{proof}

\subsection{$\mathcal B$ Graphs} Recall from Section \ref{almostplanar} that the $\mathcal B$ graphs consist of the infinite family of bicycle wheels $B_n$ and all the non-planar 3-connected deletion-minors of $B_n$. 

We begin by  proving one of the two main results in this paper, that 4-connected almost-planar graphs are pancyclic and Hamiltonian-connected. Theorem \ref{mainresult1} is repeated below for convenience.

\begin{theorem} \label{theorem1} Let $G$ be a $4$-connected  almost-planar graph.  Then $G$ is pancyclic and Hamiltonian-connected.
\end{theorem}
 
 \begin{proof} Observe that  $B_n$ is a 4-connected almost-planar graph and that $B_n\backslash z$ is 4-connected and planar. Removal of any other edge results in a graph with a degree 3 vertex, which is clearly not 4-connected.  The M\"obius ladders are cubic graphs and therefore not 4-connected. Finally, the $\mathcal H$ graphs are obtained by constructing a 3-sum of $K'''_{3,3}$ and a wheel graph, and therefore not 4-connected. Thus $B_n$ is the only 4-connected almost-planar graph. To obtain a cycle of length $k$, where $3\le k \le n$,  begin with the middle vertex $n$, then go to vertex 1 and around $k-2$ vertices along the rim, and return to vertices $n-1$ and $n$. That is, 
 $$C_k: \  n, 1, 2, \dots k-2, n-1, n.$$
In particular, $n, 1, 2, \dots, n-2, n-1, n$ is a Hamiltonian cycle. This cycle also serves as a Hamiltonian path between the two middle vertices $n$ and $n-1$, as well as a Hamiltonian path between either of the two middle vertices and a rim vertex. 
 Finally, let $i$ and $j$ be two non-adjacent rim vertices, where $i, j \in \{1, 2, \dots n-3\}$. The path 
 $$i, i+1, \dots , j-1, n, n-1, n-2, \dots , j-1, j $$ is a Hamiltonian path from $i$ to $j$. If $i$ and $j$ are adjacent, the Hamiltonian path is $$i, n, n-1, n-2, \dots , j-1, j.$$ Therefore $B_n$ is pancyclic and Hamiltonian-connected. 
 \end{proof}

Combining Theorem \ref{theorem1} with Tutte's result and Thomassen's result described in Section 1,  we may conclude that 4-connected planar and almost-planar graphs are Hamiltonian and Hamiltonian-connected.  
 
In the rest of this subsection, we will examine the non-planar  3-connected deletion-minors of $B_n$. First observe that  the edge automorphism group of $B_n$ has three orbits:  the spoke edges: the rim edges, and $z$.  Second observe  from  Figure \ref{fig2} that for $1\le i\le n-2$ and $n\ge 6$, 
$$B_n / r_{i-1} \backslash \{s_i, t_i\} \cong B_{n-1}$$
and 
$$B_n /\{r_4, \dots , r_{n-2}\} \backslash \{s_4, \dots , s_{n-2}, t_4, \dots , t_{n-2}\}\cong B_5\cong K_5.$$
Further, note that $K_5$ is minimally non-planar. Third, observe that
\begin{itemize}
\item Removing $z$  gives a planar graph; 
\item Removing any one of the rim edges $r_1, \dots , r_{n-2}$ gives a planar graph; 
\item Removing pair-wise adjacent spokes $s_i$ and $t_i$  gives a graph that is not 3-connected since the common vertex $i$ is left with degree 2;
\item Removing all the $s$-spokes (or all the $t$-spokes) gives a graph that is not 3-connected; and
\item Removing $z$ and all the $s$-spokes (or $t$-spokes), and the resulting isolated middle vertex gives a wheel with $n-2$ spokes, which is planar.
\end{itemize}

Thus the only way of obtaining a 3-connected  non-planar  deletion-minor of $B_n$, where $n\ge 6$ is to delete $s$-spokes and $t$-spokes in any combination, where $1 \le s, t \le n-4$, provided the removed spokes are not pair-wise  adjacent. 

Let us begin by examining the Hamiltonian cycles in the $\mathcal{B}$ graphs, leading to the following lemma.

\begin{lemma} \label{Ham} A $\mathcal{B}$ graph   is Hamiltonian.
\end{lemma}

\begin{proof}Since a $\mathcal{B}$ graph has at least one $s$ and one $t$ spoke, there must be two neighboring rim vertices that connect to different spokes. Let $r_i$ be the rim edge incident to spokes $s_{i-1}$ and $t_i$. A Hamiltonian cycle labeled by its edges is $$r_1 \cdots r_{i-1}s_{i-1}zt_ir_{i+1}\cdots r_{n-3}r_{n-2}.$$
\end{proof}
 
\begin{figure}[h]
\centering
\includegraphics[width=5.85in]{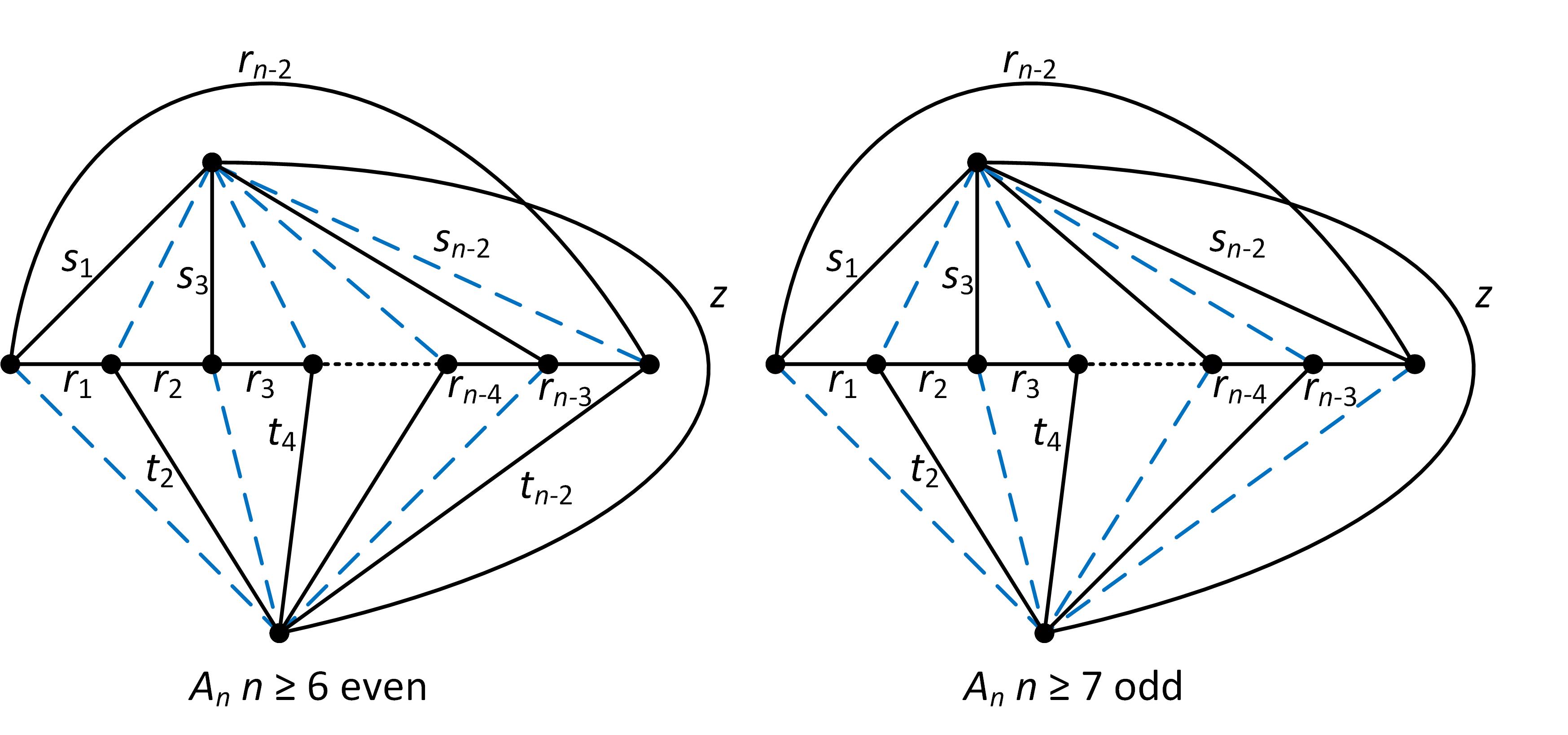}
\caption{The infinite family of graphs $A_n$, where $n\ge 6$}
 \label {An}
\end{figure}

Let $n\ge 6$ and let $1 \le i\le n-4$. Let $\mathcal G_n$ be the set of  all 3-connected non-planar deletion-minors of $B_n$ obtained by deleting $s$-spokes or $t$-spokes such that if   $s_i$ is removed, then   $t_i$ remains, and vice-versa. In other words, every rim vertex is incident to just one of $s_i$ or $t_i$. A graph in $\mathcal G_n$ is ``irreducible'' in the sense that no proper deletion-minor is both 3-connected and non-planar. Removing alternating $s$-spokes and $t$-spokes, while preserving 3-connectivity, gives the graphs displayed in Figure \ref{An}, which we will call $A_n$. 
 We will prove that the pancyclicity of $A_n$ depends on the parity of $n$.  Specifically, we will prove that if $n$ is even, then $A_n$ has only even cycles, whereas if $n$ is odd, then $A_n$ is pancyclic.

\begin{lemma} \label{Aneven}
Let $n\ge 6$ be even. Then  $A_n$ has a cycle of  length $k$ if and only if $k$ is even.
\end{lemma}

\begin{figure}[h]
\centering
\includegraphics[width=3in]{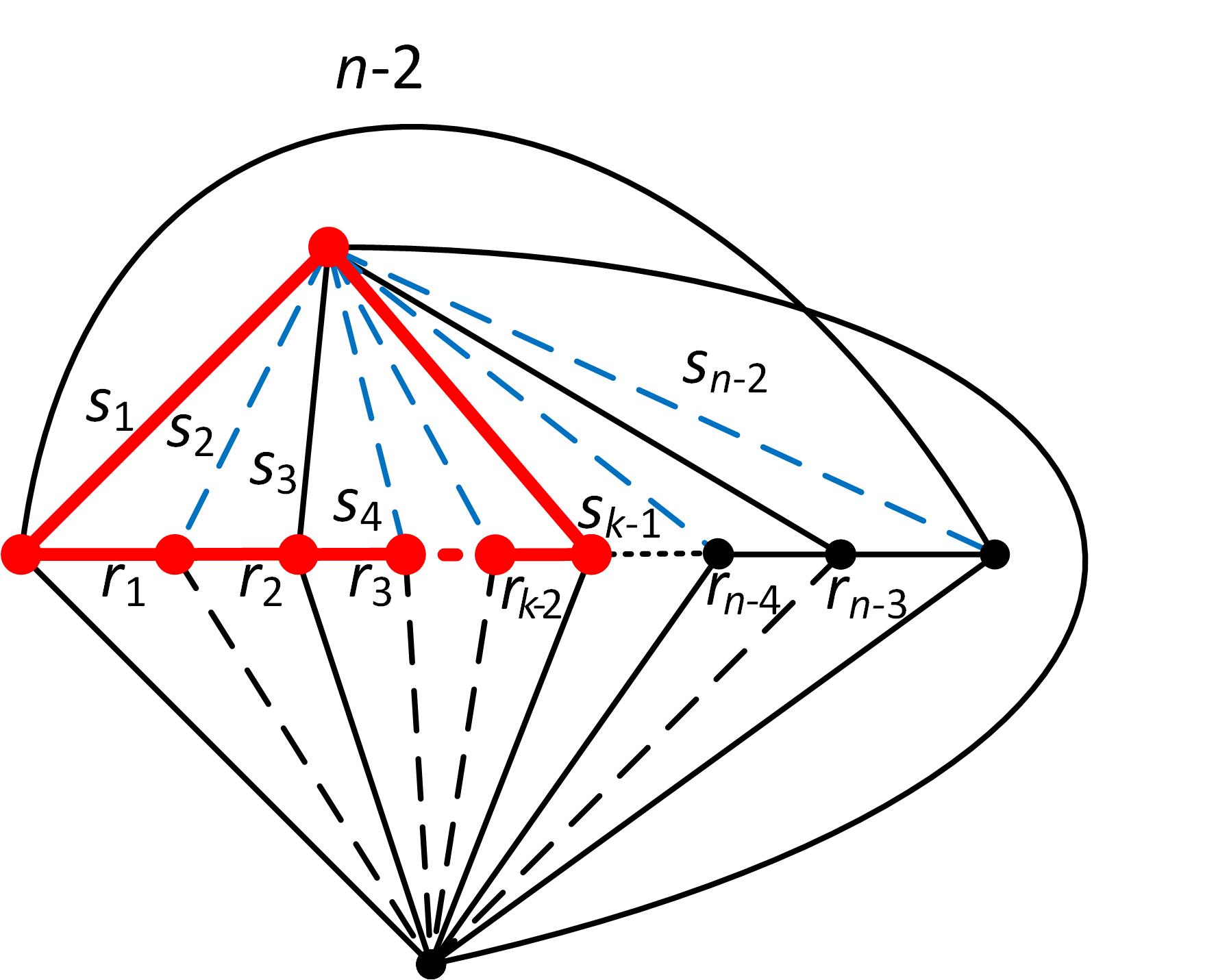}
\caption{Even cycles in $A_n$, where $n\ge 6$ is even.}
 \label {evenAnevencycles}
\end{figure}

\begin{proof} Suppose $A_n$ has a cycle of length $k$. Observe that $A_n$ is a bipartite graph since there are no edges between the two vertex classes $\{1, 3, \dots , n-3, n-1\}$ and $\{2, 4, \dots n-2, n\}$. Therefore $A_n$ has no odd cycles, and consequently $k$ is even. Conversely, suppose  $k$ is even. We will prove that $A_n$ has a cycle of length $k$. Observe from Figure \ref{An} that  the $s$-spokes labeled by even subscripts $s_2, \dots , s_{n-2}$ and the $t$-spokes labeled by odd subscripts $t_1, \dots , t_{n-3}$ are missing in $A_n$. A cycle of even length $k\in \{4, 6, 8, \dots , n-2\}$ labeled by its edges is shown in Figure \ref{evenAnevencycles} and obtained as follows: $$s_1r_1 r_2   \cdots  r_{k-2} s_{k-1}.$$
Lastly, a Hamiltonian cycle labeled by its edges is
$$s_1r_1t_2t_{n-2}r_{n-2}\cdots r_3s_3.$$
\end{proof}

The next result is the main result of this subsection. Observe that a $\mathcal B$ graph with $n$ vertices is obtained by adding spokes to graphs in $\mathcal G_n$. We will prove that a graph in $\mathcal G_n$ with two adjacent $s$-spokes (or two adjacent $t$-spokes) is pancyclic. This argument will cover $A_n$, where $n$ is odd, since spokes $s_1$ and $s_{n-2}$ are adjacent. What will be left to analyze is $A_n$, where $n$ is even, and we already showed in Lemma \ref{Aneven} that it is not pancyclic. It will then follow that every $\mathcal B$ graph with the exception of $A_n$, where $n$ is even, is pancyclic. 

\begin{theorem} \label{Bgraphs}
Let $n\ge 5$ and let $G$ be a $\mathcal B$ graph on $n$ vertices such that $G\not\cong A_n$, where $n$ is even. Then $G$ is pancyclic.
\end{theorem}

\begin{figure}[h]
\centering
\includegraphics[width=5.25in]{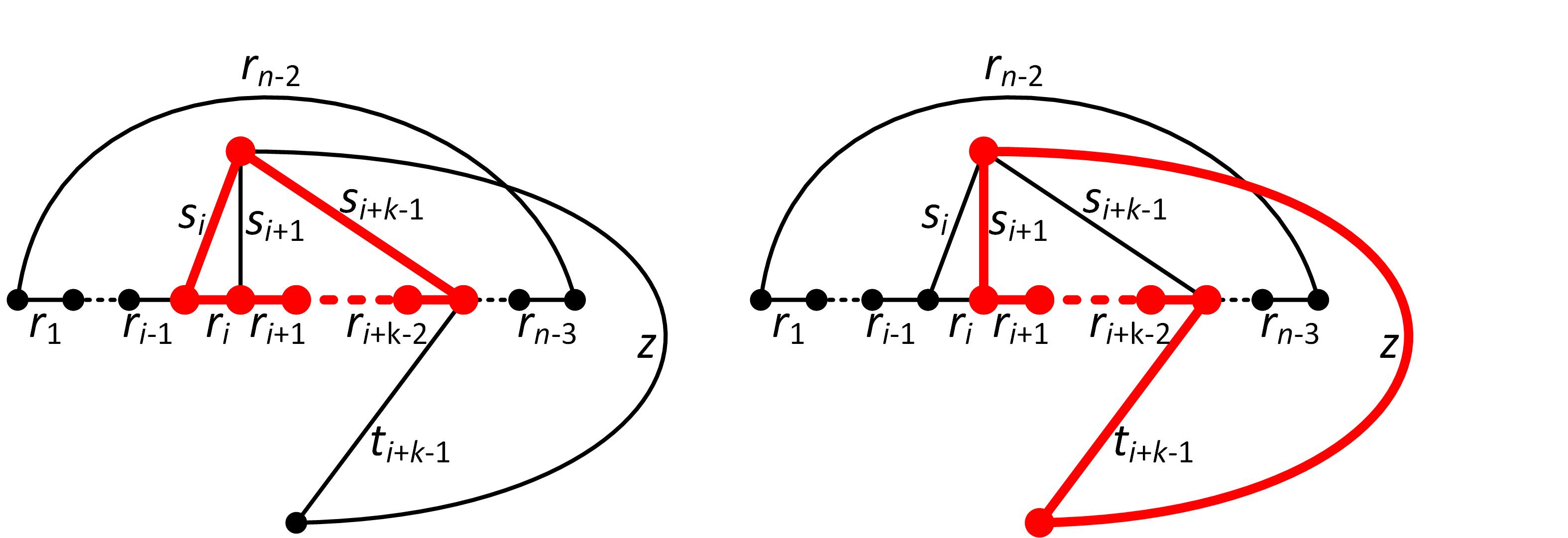}
\caption{Cycles when two adjacent spokes are present.}
 \label {twoadjspokes}
\end{figure}

\begin{proof}   
Lemma \ref{Ham} implies that a $\mathcal B$ graph has a cycle of size $n$. So we may focus on cycles of size strictly less than $n$. Since $G$ is a $\mathcal B$ graph on $n$ vertices, $G$ is obtained from a graph in $\mathcal G_n$ by adding one or more of the missing spokes. It suffices to examine the pancyclicity of the graphs in $\mathcal G_n$. We will prove that the only graphs in $\mathcal G_n$, that are not pancyclic are $A_n$, where $n$ is even.  

Let $G\in \mathcal G_n$ and suppose that $G$ has two adjacent $s$-spokes or two adjacent $t$-spokes. Without loss of generality suppose $s_i$ and $s_{i+1}$ are the adjacent spokes, where $1 \le i \le n-2$ and the subscripts are taken modulo $n-2$, that is $s_{n-2+1}=s_1$. Then $s_ir_is_{i+1}$ is a 3-cycle. 
Next, let $4\le k \le n-1$. If spoke  $s_{i+k-1}$ is present in $G$, then as shown in Figure \ref{twoadjspokes}, 
$$s_ir_ir_{i+1} \cdots r_{i+k-2}s_{i+k-1}$$ 
is a $k$-cycle.
If spoke  $s_{i+k-1}$ is not present in $G$, then spoke $t_{i+k-1}$ is necessarily present, and in this case, 
$$s_{i+1}r_{i+1} \cdots r_{i+k-2}t_{i+k-1}$$
is a $k$-cycle.

Thus we may assume that $G$ does not have two adjacent $s$ spokes or two adjacent $t$ spokes. Since $G$ is 3-connected and non-planar, $G \cong A_n$. If $n\ge 5$ is odd, then $A_n$ has two adjacent spokes and is therefore pancyclic. So we may assume $G\cong A_n$, where $n$ is even. Lemma \ref{Aneven} implies that $A_n$ has a cycle of  length $k$ if and only if $k$ is even. 

To complete the proof observe that adding any missing spoke to $A_n$, where $n$ is even, gives a graph with two adjacent spokes, which is pancyclic.  
\end{proof}
 
\subsection{$\mathcal H$ Graphs} 

 For $p, q, r\ge 1$,  ${\mathcal H_1}(p, q, r)$ and ${\mathcal H_2}(p, q, r)$ are constructed from $K'''_{3,3}$ by attaching fans along triangles in the specific ways shown in Figure \ref{fig4}. When the context is clear we will refer to these families as ${\mathcal H_1}$ and ${\mathcal H_2}$.

A 3-connected graph is \textit{minimally $3$-connected} if removal of any edge destroys 3-connectivity. Observe that most of the edges in $\mathcal H_1$ and $\mathcal H_2$ are incident to degree 3 vertices and therefore cannot be deleted. In both families, only edges $ab$, $bc$ and $ac$ may be deleted. We will prove that besides $K_{3,3}$, the minimally 3-connected infinite families of graphs  ${\mathcal H_1}\backslash \{ab, bc, ac\}$ and 
${\mathcal H_2}\backslash \{ab, bc, ac\}$ are pancyclic and therefore so are ${\mathcal H_1}$ and ${\mathcal H_2}$.

\begin{figure}[h]
\centering
\includegraphics[width=6.8in]{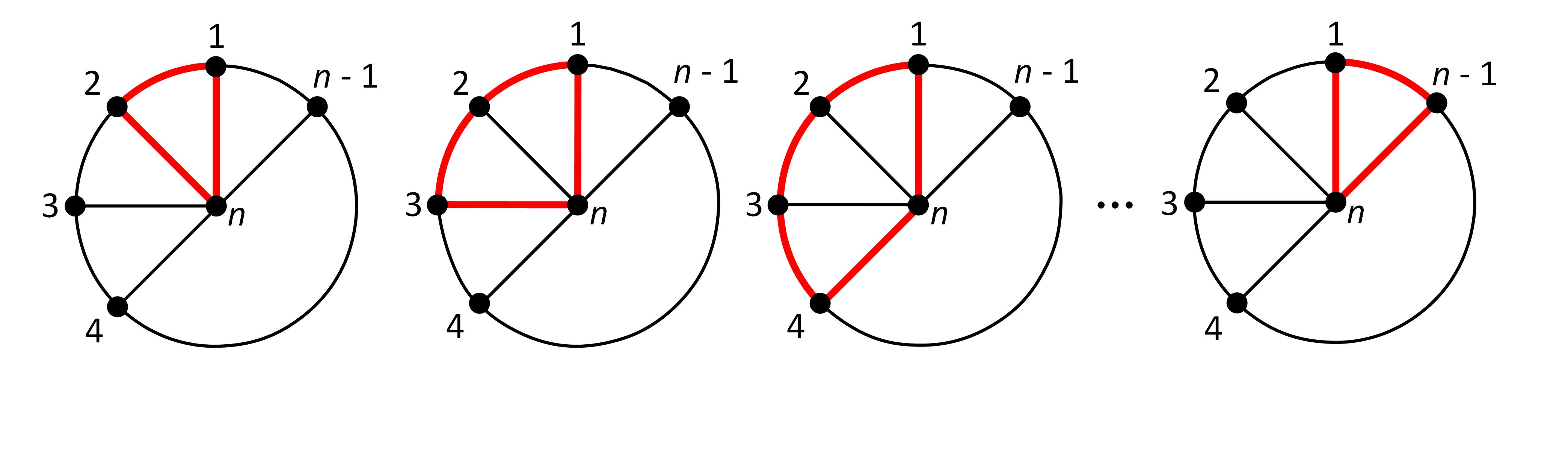}
\caption{The wheel graph is pancyclic.}
 \label {Wn}
\end{figure}

A  key idea  in the proof is the pancyclicity of the wheel. Cycles of all sizes in $W_{n-1}$ are illustrated in Figure \ref{Wn}, where a cycle of length $t$ is obtained by either selecting some spokes or avoiding some spokes.
For example, a cycle of length 3 is obtained by selecting 2 spokes or avoiding $(n-1)-2=n-3$ spokes.  In general, for $3\le t \le n-1$, a cycle of length $t$ is obtained by selecting $t-1$ spokes or avoiding $n-1-(t-1)=n-t$ spokes.  Figure \ref{Wn} and the following table make this idea clearer. We use both approaches in next proof.

\begin{table}[h]
\centering
\begin{tabular}{|c|c|}
\hline
$C_3: 1, n, 2 1$  & $C_n:  1,  n, 2, 3, \dots , n-1, 1$ \\ 
\hline
$C_4: 1, n, 3, 2 1$ & $C_{n-1}: 1, n, 3, 4, \dots , n-1, 1$ \\  
\hline
$C_5: 1, n, 4, 3, 2, 1$ & $C_{n-2}: 1, n, 4,   \dots , n-1, 1$ \\ 
\hline  
$\vdots$ & $\vdots$ \\  
\hline  
$C_n: 1, n, n-1 \dots , 3, 2, 1$ & $C_3: {\bf 1, n, n-1}, 1$ \\
\hline
\end{tabular}
\end{table}

Another idea used in the proof is as follows: consider an arbitrary graph $G$ with a triangle $e, f, g$, and suppose we attach a fan along sides $e$ and $f$ as shown in Figure \ref{fig3}. Further suppose $G$ is pancyclic and $G$ has a Hamiltonian cycle $C$ that uses one of $e$ or $f$. Then the new graph obtained from $G$ by attaching a fan is also pancyclic.   Looking at cycles in this manner provides intuition for why under certain conditions, a graph with a Hamiltonian cycle is also pancyclic, even when most of the vertex degrees are small.

\begin{theorem} \label{Hgraphs}
Let $G$ be an $\mathcal H$-graph such that $G\not \cong K_{3,3}$.  Then $G$ is pancyclic.
\end{theorem}

\begin{proof} First, observe that $K_{3,3}$ is bipartite, and therefore has no odd cycles. 
However, $K_{3,3}'$ is pancyclic. To see this, view $K_{3,3}'$ as the graph obtained from $K'''_{3,3}$ displayed as the first diagram of Figure \ref{fig4} with edges $bc$ and $ac$ removed. Cycles of lengths 3, 4, 5, and 6 are:
$$C_3: a, x_1, b, a$$ $$C_4: a, x_1, b, y_1, a$$ $$C_5: a, b, x_1, c, y_1, a$$ $$C_6: a, x_1, b, y_1, c, z_1, a.$$  Consequently, adding additional edges to obtain  $K_{3,3}''$ and $K_{3,3}'''$  preserves  pancyclicity. We will prove that $K_{3,3}$ is the only $\mathcal H$ graph that is not pancyclic. Note that the number of vertices in $\mathcal H$ graphs is $n=p+q+r+3$.

Call a spoke a $p$-spoke if its end vertices are $bx_i$ for $1 \le i \le p$; a $q$-spoke if its end vertices are $by_j$ for $1 \le j \le q$; and a $r$-spoke if its end vertices are $az_k$ for $1 \le k \le r$.

\textit {Claim 1:  $\mathcal{H}_1\backslash \{ab, \ bc, \    ac\}$, where at least one of $p, q, r$ is greater than $1$, is pancyclic.}

The infinite family ${\mathcal H_1}$,  shown in Figure \ref{fig4}, is obtained by attaching fans of length $p$, $q$, and $r$ in triangles $\{a, b, x_1, a\}$,  $\{b, c, y_1, b\}$,   and  $\{a, b,  z_1, a\}$ of $K'''_{3,3}$ along sides $ab$ and $bx_1$;  $bc$ and  $by_1$; and $ab$ and $bz_1$,  respectively.  Observe from Figure \ref{H1} that 
 $$x_1,  \dots , x_p, a, z_r, \dots z_1, c, y_q, \dots y_1, b, x_1$$ 
is a Hamiltonian cycle of   $\mathcal{H}_1 \backslash \{ab, \ bc, \    ac\}$.

\begin{figure}[h]
\centering
\includegraphics[width=2.15in]{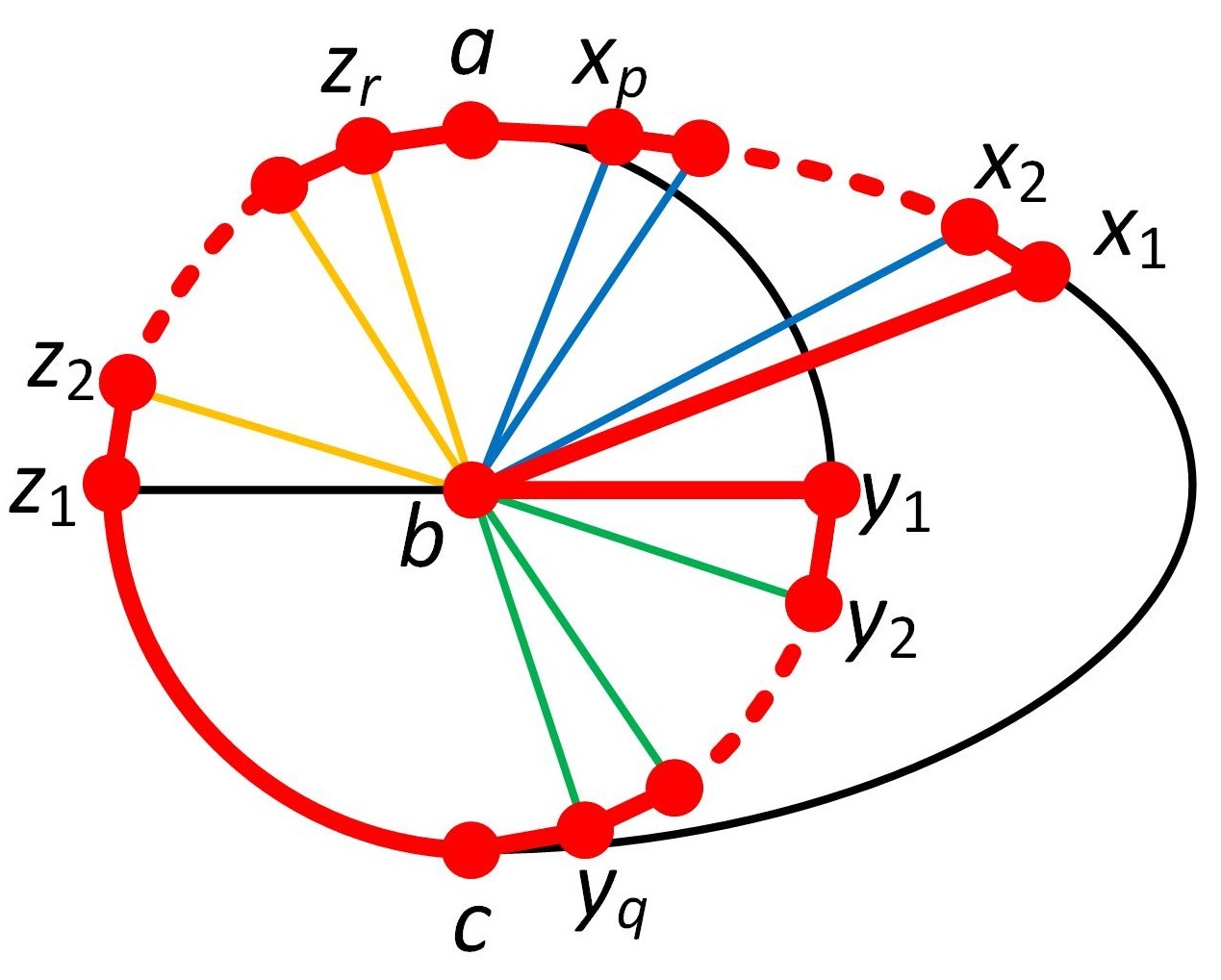}
\caption{Pancyclicity of $\mathcal{H}_1(p,q,r)$}
 \label {H1}
\end{figure} 

Moreover, the deletion-minor $\mathcal{H}_1\backslash \{ab, \ bc, \  ac\} \backslash \{ ay_1, \ cx_1\}$  is isomorphic to the wheel graph on $n$ vertices  without one rim edge $x_1y_1$ and two spokes $ab$ and $bc$. Following the procedure described in Figure \ref{Wn}, we can obtain cycles of all but six sizes. 
Since spoke $ab$ is missing, this procedure will not give a cycle of size $p+2$ and $n-(p+2)=p+q+r+3-(p+2)=q+r+1,$ and $r+2$ and $n-(r+2)=p+q+r+3-(r+2)=p+q+1$.
Since spoke  $bc$ is missing, this procedure will not give a cycle of size  $q+2$ and  $n-(q+2)=p+q+r+3-(q+2)=p+r+1$.  The  missing cycle lengths are obtained in a slightly different manner as indicated in the following table:

\begin{table}[h]
\centering
\begin{tabular}{|c|c|}
\hline
{\bf Cycle Length}  & {\bf Cycle} \\ 
\hline
$p+2$  & $b, x_2, x_3, \dots , x_p, a, z_r, b $ \\ 
\hline
$q+2$   & $b, y_2, y_3, \dots , y_q,  c, z_1, b$ \\ 
\hline
$r+2$   & $b, z_2, z_3, \dots , z_r,  a, x_p, b$    \\ 
\hline
$q+r+1$   &  $b, y_2, y_3, \dots , y_q, c, z_1, z_2, \dots , z_r, b$  \\ 
\hline
$p+q+1$   &  $b, x_1, x_2, \dots , x_p, a, y_1, y_2, \dots , y_q, b$   \\ 
\hline
$p+r+1$   &   $b, x_2, x_3, \dots , x_p, a, z_r, \dots , z_1, b$ \\ 
\hline
\end{tabular}
\end{table}

This completes the proof of Claim 1.

\textit {Claim 2:  $\mathcal{H}_2\backslash \{ab, \ bc, \    ac\}$, where at least one of $p, q, r$ is greater than $1$, is pancyclic.}

The infinite family $\mathcal{H}_2$,  shown in Figure \ref{fig3},  is obtained by attaching fans of length $p$, $q$, and $r$ in triangles $\{a, b, x_1, a\}$,  $\{b, c, y_1, b\}$,   and  $\{a, b,  z_1, a\}$    along sides $ab$ and $bx_1$;  $bc$ and  $by_1$; and $ab$ and $az_1$,  respectively.   Observe from Figure \ref{H2} that 
$$a, x_p, \dots , x_1, b, y_1, \dots , y_q, c, z_1 \dots , z_r, a.$$  
is a Hamiltonian cycle of   $\mathcal{H}_2 \backslash \{ab, \ bc, \    ac\}$.
This cycle  goes along one side of each of three fans. In particular $C$ goes through $bx_1$, $by_1$, and $az_r$. 

\begin{figure}[h]
\centering
\includegraphics[width=2.15in]{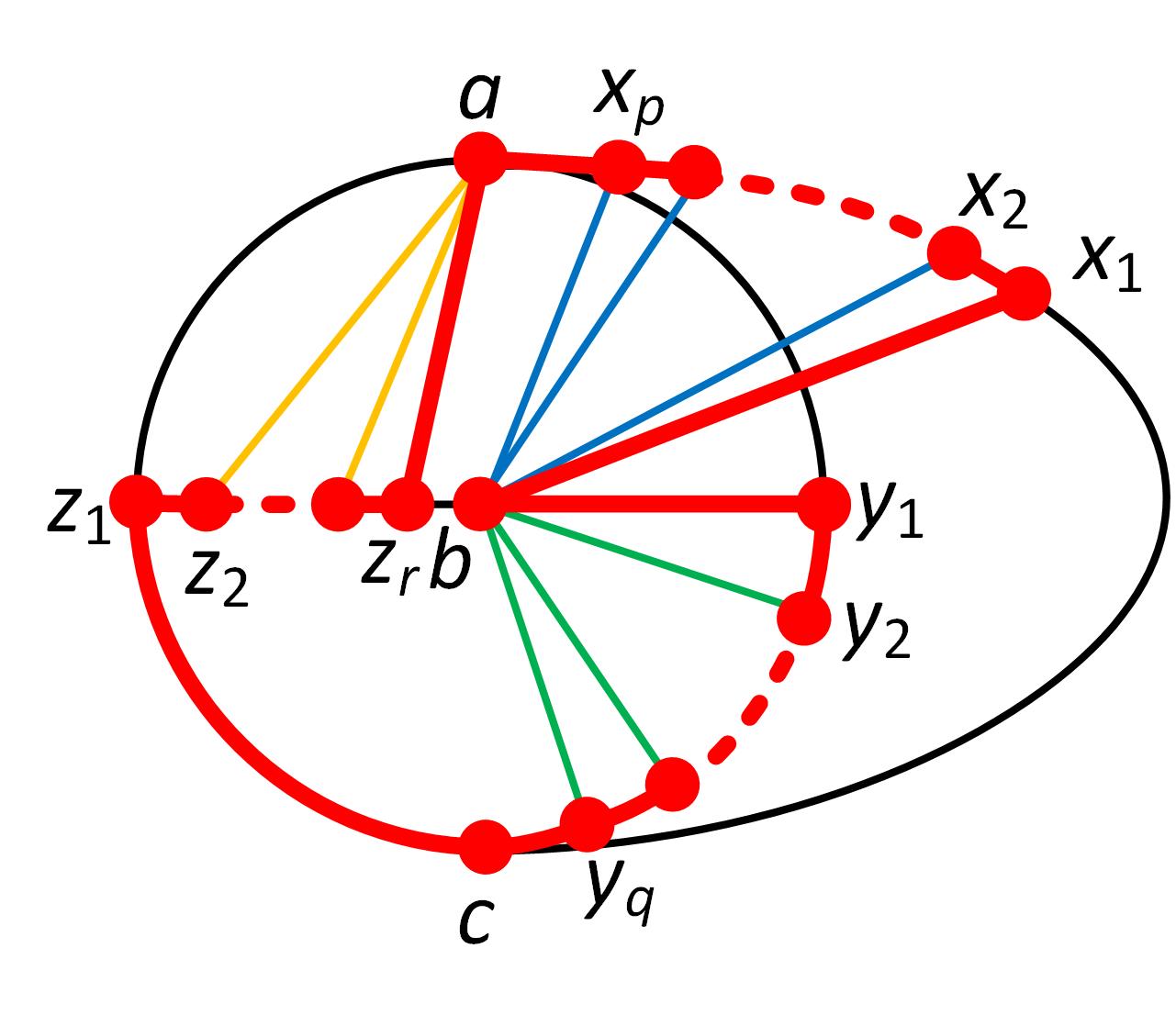}
\caption{Pancyclicity of $\mathcal{H}_2(p,q,r)$}
 \label {H2}
\end{figure}

A cycle of length $(p-i)+q+r+3$ may be obtained by avoiding $i$ consecutive $p$-spokes $bx_1, \dots bx_i$, where $1\le i\le p-1$, since at least one $p$-spoke must remain to preserve edge $ax_p$. This cycle of length $n-i$ is  $$a, z_r, \dots, z_1, c, y_q, \dots , y_1,  b, x_{i+1}, \dots, x_p, a.$$
A cycle of length $p+(q-j)+r+3$ may be obtained by avoiding $j$ consecutive $q$-spokes $by_1, \dots by_j$, where $1\le j\le q-1$, since at least one $q$-spoke must remain to preserve edge $cy_q$.  
This cycle of length $n-j$ is   $$a, z_r, \dots, z_1, c, y_q, \dots ,   y_{j-1},  b, x_1, \dots, x_p, a.$$
A cycle of length $p+q+(r-k)+3$ may be obtained by avoiding $k$ consecutive $r$-spokes $az_2, \dots az_k$, where $1\le k\le r-1$ since at least one $r$-spoke must remain. In this case when $r-1$ spokes are removed edge $ar_1$ remains.  
This cycle of length $n-k$ is  $$a, z_{k+1}, \dots, z_1, c, y_q, \dots , y_1, b, x_1, \dots, x_p, a$$
Since the Hamiltonian cycle goes along one side of each of the three fans, any combination of $i$, $j$, and $k$ vertices may be avoided to obtain a cycle of length $n-(i+j+k)$.  This completes the proof of Claim 2.

Finally, $\mathcal{H}_1 \backslash \{ab, bc, ac\}$ and $\mathcal{H}_2 \backslash \{ab, bc, ac\}$ are pancyclic, so are $\mathcal{H}_1$ and  $\mathcal{H}_2$.  Therefore, $K_{3,3}$ is the only $\mathcal{H}$-graph that is not pancyclic.
\end{proof}

Now, we put all the pieces in this subsection together to obtain the following result:

\begin{theorem}\label{mainresult3}
Let $G$ be a $3$-connected almost-planar graph such that $G\not \cong  A_n$ for even $n\ge 6$  and $G\not \cong V_{2k}$ for $k\ge 4$. Then $G$ is pancyclic.
\end{theorem}

\begin{proof}
Theorem \ref{Mobiusgraphs} implies that $V_{2k}$ for $k \ge 4$ is not pancyclic.
 Theorem \ref{Bgraphs} implies that all $\mathcal B$-graphs are pancyclic except  $A_n$ for even $n\ge 6$. Theorem \ref{Hgraphs} implies that all $\mathcal H$-graphs are pancyclic except   $K_{3,3}$.  The result follows since $A_6\cong K_{3,3}$.
\end{proof}

We are now ready to prove Theorem \ref{mainresult2}, which we state again for convenience.

\begin{theorem}\label{mainresult3}
A $3$-connected almost-planar graph  is pancyclic if and only if it has a cycle of length $3$.
\end{theorem}

\begin{proof} Suppose $G$ is a $3$-connected almost-planar graph that is pancyclic.  Then clearly $G$ has a cycle of length 3. The converse is what we have to prove; that is, if $G$ is a $3$-connected almost-planar graph with a $3$-cycle, then $G$ is pancyclic. Since $G$ is a 3-connected almost-planar graph, by Theorem \ref{almostplanartheorem} $G\cong V_{2k}$, for $k\ge 4$, or $G$ is a $3$-connected non-planar deletion-minor of $B_{n}$,  for $n\ge 5$, or ${\mathcal H_1}(p, q, r)$ or  ${\mathcal H_2}(p, q, r)$, for $p, q, r \ge 1$. 

Theorem \ref{Mobiusgraphs} implies that if $k\ge 4$ is odd, then the M\"obius ladders $V_{2k}$ have only cycles of even length, whereas if $k\ge 4$ is even, then $V_{2k}$ have cycles of even length as well as, cycles of length greater than $k$. In particular, the smallest M\"obius ladder $V_8$ has cycles of length 4, 5, 6, 7, and 8. There is no M\"obius ladder with a cycle of size 3, and this is the only cycle length missing among all the M\"obius ladders.

Next, consider the infinite family of irreducible graphs $A_n$, where $n\ge 5$. If $n$ is even, then $A_n$ is bipartite and therefore has only  cycles of even length. In particular, it has no cycle of length 3. Adding any spoke to $A_n$ immediately gives a 3-cycle.   If $n$ is odd, then  $A_n$ has a 3-cycle as shown in Figure \ref{An}.  Theorem \ref{Bgraphs} implies that the only $\mathcal B$ graph that is not pancyclic is  $A_n$, for even $n\ge 6$, and that all the other $\mathcal B$-graphs are pancyclic. 

Finally, consider the $\mathcal H$ graphs. Observe that $K_{3,3}$ has only even cycles. All other $\mathcal H$-graphs have a 3-cycle, and by Theorem \ref{Hgraphs} they are pancyclic.  
\end{proof}

\section{\bf Conclusion}

The class of non-planar graphs contains the bicycle wheels and their non-planar deletion-minors. These graphs form a very large and important class of 3-connected graphs in graph structure theory, and we proved that most of them are pancyclic. The next step is to analyze the planar deletion-minors of bicycle wheels. This appears to be a more difficult, but still tractable problem. In a similar manner, we can continue examining minimally 3-connected graphs for Hamiltonicity by finding other large classes with many low degree vertices that are Hamiltonian, and then examining them for pancyclicity. Such graphs are not covered by prevailing degree-based pancyclicity conditions.

The second author has developed structural results for generating cyclically 4-connected graphs, especially cyclically 4-connected planar graphs. The next step in this direction is to take each one of these graphs and construct their line graphs, thereby obtaining a large number of planar graphs with girth at least 5 to see if they are missing any other cycles. 

Third, we already know based on existing results that a 4-connected planar graph with $n$ vertices has large cycles, that is cycles of size at least $n-6$. Based on a preliminary investigation, we believe existing techniques can be extended to show that they will also contain at least one smaller size cycle.

A fourth direction is to develop structural results for 4-connected planar graphs in order to generate and examine them for pancyclicity.


\bigskip
{\bf Acknowledgements:} This research was partially supported by a National Science Foundation grant (NSF DMS 2051026) and a Caltech SURF award.

{\bf Declaration} The authors have no relevant financial or non-financial interests to disclose.


\end{document}